\documentclass[11pt, twocolumn]{article}
\usepackage{article}
\author{Torgeir Aamb\o\thanks{Norwegian Defence Research Establishment (FFI)}}
\title{Collective deterrence as a classification problem: \\
Voting rules, deterrence credibility, and escalation risk}
\date{}

\begin{document}
\maketitle

\begin{abstract}
    \textbf{Deterrence coalitions that collectively own their deterrence technology, need an institutional design to decide when to retaliate against an attack or incident. This choice of institutional design, formalized through a social choice function, introduces a tradeoff between credible deterrence and escalation risk. We study this tradeoff via a simple signalling model, and use it to construct an associated binary classification problem to determine institutional designs that perform well in a variety of environments. For a small coalition of four members, we compute and study the statistics of the empirical ROC curves associated to a variety of choice functions and probability distributions for retaliation and false positives.} 
\end{abstract}

\emph{Convince your enemy that he will gain very little by attacking you, which will diminish their enthusiasm.} \\ 
\hfill-- Sun Tzu, \emph{The Art of War} \cite{Art_of_war}.

\section{Introduction}
Since the emergence of modern warfare, understanding how to deter potential attackers has become a central concern of strategic studies and a core element of state security. Deterrence as a military strategy -- avoiding an attack by maintaining sufficiently powerful means of guaranteed retaliation \cite{Schelling_1966,Snyder_1961} -- rests on shaping an adversary's beliefs such that the expected costs outweigh the perceived benefits.

In this simplified form there are two strategies, reduce the potential benefits or increase the costs. The former is often counter-productive, though has been pursued through tactics such as scorched-earth, see for example \cite{Kreike_2021}. The latter, however, requires that the potential attacker believes the consequences to be large, for example by believing that large scale retaliation would occur. 

In essence, deterrence relies on a push-pull between two fundamental objectives \cite{Schelling_1960}: 
\begin{enumerate}
    \item the resolve to use retaliatory means has to be sufficiently credible; 
    \item one should limit the possibility and risk of escalation spirals and accidental use. 
\end{enumerate} 
The dynamic between these two mechanisms is especially potent for nuclear deterrence, where the consequences are potentially of an existential nature. For a classical outlook on these fundamental issues see \cite{Schelling_1960} and \cite{Snyder_1961}. Classical deterrence theory does, however, not often examine coalitions and their collective deterrence capabilities, which is the goal of this paper. 

Currently, the only proper nuclear alliance is NATO, in the sense of being a multinational alliance and not a bilateral agreements, where the deterrence capabilities of individual nations are also shared to allies. However, NATO does not own nuclear weapons as a collective entity, nor does it operate a jointly owned nuclear command structure with shared launch authority -- the historical proposals for a NATO multinational nuclear force never materialized \cite{Trachtenberg_1999}. For example, U.S. weapons stationed in Europe remain under exclusive U.S. custody and control, and any nuclear mission requires authorization from the relevant national leaders in addition to political approval through NATO's Nuclear Planning Group (NPG). 

Hence, there is currently no coalition of states in the world that jointly own both the nuclear deterrence capabilities and the command structure to act upon threats. However, such coalitions might become more politically viable now that arms-control frameworks, such as the New START treaty \cite{NewSTART_2010}, are intensifying the European discussions about burden sharing and the possible future role of France or the UK in European nuclear deterrence. 

Currently, of course, the Non-Proliferation Treaty (\cite{NPT_1968}) obstructs the creation of such new nuclear alliances. However, there is room to break the treaty in the case of existential threats, or, alternatively, one could consider other deterrence technologies not limited by such restrictions, see \cite{EuropeanNuclearPlanningGroup}. In the advent of such coalition then, the additional problem of \emph{collective} rather than \emph{unilateral} decision-making arises to complicate the above two problems. The problem of \emph{credible resolve} becomes a problem of how the coalition's states aggregate their individual resolves into a joint decision, and the problem of \emph{escalation risk} becomes a problem of not allowing this joint decision to be too easily affected by individual potential errors. Hence this joint institutional decision-structure introduces a trade-off between the two fundamental problems. 

We can then phrase the motivating question for the present paper as: How does the institutional design of a coalition affect the deterrence mechanism and its effectiveness? In this paper we study this question as a classification problem using Receiver Operating Characteristic (ROC) curve analysis.

Aggregations of individual choices are in social choice theory measured by a \emph{social choice function}, see for example \cite{Arrow_1951, Moulin_1988}. To set up our method, we develop a simple Bayesian decision model inspired by \emph{signalling games} \cite[Part IV]{Fudenberg-Tirole_1991}. This model is dependent on a weighted threshold choice function $f$ that aggregates the coalition's votes on retaliation into a shared resolve. 
 
To understand the interactions between deterrence and institutional design, we compute simple attack conditions, and use this to formulate the choice of institutional design for the coalition as a \emph{binary classification problem}. Each institutional design has an associated empirical ROC curve -- see \cite{Fawcett_2006} -- for each information environment, formalized by a list of probabilities for resolves and false positives for the coalition members. By studying the statistics of these ROC curves for a small hypothetical coalition of four members, we can study which designs that perform well and which do not. This does not capture the full strategic payoffs for the deterrence coalition, but gives a tractable and computable way to compare social choice functions on their true positive-false positive detection abilities. 

A voting scheme that performs well on a variety of probability distributions for retaliation and false alarms is, perhaps not surprisingly, the unbiased symmetric one, where each coalition members has equal say. Close behind follows a scheme where each coalition members is given a weighted vote according to some relevant capabilities, like ISR. Across our stylized information environments, the worst performer is the dictatorial scheme, where only one country has voting power over the use of deterrence technology.

Finally we use Youden's $J$-statistic (\cite{Youden_1950}) to study which thresholds that give the best tradeoff between credible deterrence and escalation risk, concluding that the best threshold usually is the \emph{majority} threshold for our small hypothetical coalition of four members. 

Our goal with this paper is to provide a first step toward understanding the interaction between institutional design and collective deterrence by introducing, studying and exploring a novel binary classification-based framework. We deliberately adopt a simplified game model to isolate the informational tradeoffs inherent in collective decision-making -- extensions incorporating more strategic behavior, correlated information, and proofs of optimality will be studied in future work. 

\subsection{Contributions}
This paper has three main contributions towards the understanding of institutional design in collective deterrence: 
\begin{enumerate}
    \item introducing a novel classification-based framework for analyzing collective deterrence, with an easily understandable visual element;
    \item providing a tractable model for linking institutional design to deterrence credibility and escalation risk;
    \item illustrating the framework through simulations of small coalitions under varying informational environments.
\end{enumerate}

Our results should be interpreted as illustrative patterns and not general results. The results, however, give clear directions for future work on optimal institutional designs using a more general and proof-based approach.

\section{Theoretical background}
Our model is based on a combination of two much studied mathematical models: signalling games and social choice functions. These are each ubiquitous in their own respective fields of deterrence theory and voting theory. Before putting forth our combined model, we briefly explain the pieces saparately.

\subsection{Signalling games}
Signalling games are dynamic Bayesian games, first studied in relation to learning by Lewis \cite{Lewis_1969}. For a more modern treatment see \cite{Fudenberg-Tirole_1991}, and for deterrence specific signalling theory see \cite{Schelling_1960}. 

In a signalling game, the players have incomplete information about the initial situation, and the game unfolds sequentially over time, with actions known to the other players \cite{Fudenberg-Tirole_1991}. More precisely, one player -- often called the sender -- has some private information about the state of the world, and takes an observable action $S$, called the signal, after which another player -- often called the receiver -- updates their beliefs about the world and chooses an action. The private information of each player is called the \emph{type}. 

In our case of deterrence, incomplete information is essential; the uncertainties around the opponents resolve lies at the heart of the mechanism making deterrence a viable strategy. 

In non-cooperative game theory, Nash equilibria are the standard solution concept for determining optimal strategies. In signalling games, one usually needs more refined concept, due to the players not knowing the beliefs of the other players correctly. These are called \emph{perfect Bayesian equilibria}, see \cite{Fudenberg-Tirole_1991_paper}. We will, however, not solve for this solution concept in this paper, as our game will be of a deliberately simple nature where this full generality is unnecessary. 

Instead we will use \emph{attack conditions}, which are assignments of simple actions to the type of the players. These are motivated by threshold logic in crisis bargaining under incomplete information, see for example \cite{Fearon_1995}, and is a simpler way of studying basic strategic actions. 

\subsection{Social choice functions}
\label{sec:social-choice-functions}

Social choice functions are mathematical models for ranking a set of alternatives, often called \emph{social states}, by their collective desirability -- combining the independent preferences of the individual into a common ``will'' of the people. Early versions were studied by Condorcet and Llull, see \cite{Condorcet_1785} and \cite{Colomer_2013}, but the modern theory now in use is often attributed to Arrow's work \cite{Arrow_1951} -- see \cite{Moulin_1988} for an introduction. 

Social choice functions can be either \emph{cardinal} or \emph{ordinal} \cite{Sen_1977}. Ordinal functions only use the relative ordering of social states to compute the best preferred state according to some mathamatical rule. On the other hand, cardinal functions use specific values attached to each social state, and hence is able to compute not just which state is preferred, but how much it is preferred to the other states. It is the former, the ordinal choice functions, that are most often used in public voting and electoral systems, and also the type of function we will use in this paper. Our needs are also very simple, as our list of potential outcomes will be binary, hence we focus only on these in the following examples. 

The most frequently used voting rule is the \emph{majority rule}, stating that the preferred choice between two options is the one chosen by more than half of the voters. By May's Theorem \cite{May_1952}, this is the only \emph{fair} voting scheme. In mathematical terms, given a list of votes $v$, we have:
\[f(v) = 
\begin{cases}
    1, &\displaystyle \sum_{i=1}^{n} v_i > \frac{n}{2}  \\
    0, & \text{else}.
\end{cases}
\]
By varying the threshold for the rule outputing a $1$, we also vary the choice function. For example, the consensus rule, used for example when NATO votes on shared policies and solutions \cite{NATOConsensus}, states that the preferred choice between two options is the one chosen by \emph{all} voters:
\[f(v) = 
\begin{cases}
    1, &\displaystyle \sum_{i=1}^{n} v_i = n  \\
    0, & \text{else}.
\end{cases}
\]
This is essentially a voting scheme in which every voter has a veto right. Choosing other thresholds, like $\frac{2}{3}$ gives the \emph{qualified majority rule} and choosing $\frac{3}{4}$ gives the \emph{supermajority rule}. 

There are also many other voting schemes, like the \emph{dictatorial} scheme, where only one voter has a say and always completely determines the outcome. The restricted class of social choice functions used in the model is described in \cref{ssec:weighted-threshold-functions}. 

\section{The model}
\label{sec:the-model}
Our model adopts this basic structure of a signalling game, as presented above, where we assume one of the players to be a coalition consisting of multiple states $\C = \{1, 2, \dots, n\}$. The other player is some potential attacking adversarial actor $\A$. Before we explain the individual parts, we summary the sequence of events: 
\begin{enumerate}
    \item the coalition has a predetermined fixed social choice function $f$, which is publicly known;
    \item nature draws aversary type $T$;
    \item nature draws signals $s_i$ for each coalition member;
    \item adversary actor $\A$ chooses to attack or not based on the above information;
    \item in the case of an attack, the coalition votes on retaliation using $f$. 
\end{enumerate}

\subsection{Types and signals}
The adversary $\A$ has a \emph{type}, $T \in \{0,1\}$, where $1$ represents being aggressive -- actively planning an attack -- and $0$ represents a benign adversary with no such intentions. The coalition has a shared prior belief $\pi = \Pr(T=1)$ of which type $\A$ has. 

Each member $i$ of the coalition $\C$ has a binary signal $s_i \in \{0,1\}$, which is understood as a reduced-form decision-relevant indicator, summarizing both informational assessments about the adversary $\A$ and the member's willingness to retaliate. This is dependent on their belief of the type of $\A$, giving a probability
\[p_i = \Pr(s_i = 1 \mid T=1)\]
for whether coalition member $i$ will vote to retaliate or not. Similarly, coalition members might receive false alarms, which gives a probability
\[q_i = \Pr(s_i = 1 \mid T=0)\]
of a retaliatory strike based on erroneous information. Note that an adversary $\A$ of type $T=0$ does not attack, so any such retaliation is an unwarranted escalation. 

\begin{remark}
    We assume here that signals and probabilities are independent across members, even though this is highly unlikely in a real scenario. In an actual coalition, information and intelligence is shared, discussions are held, and decision-makers are influenced by other members. This said, the choice allows us to gain computational power and simplicity. With our focus on understanding the underlying \emph{mechanisms} and not be as true-to-reality as possible, we feel that this simplification is sufficiently justified. 
\end{remark} 

In the case of an attack, each coalition member $i\in \C$ then decides on a binary vote $v_i \in \{0,1\}$, interpreted as their vote towards authorizing retaliation. In order to focus this paper on the institutional mechanisms, and not on whether the coalition members try to strategize for their own gain, we assume that $v_i = s_i$, which means that we assume the actors to vote truthfully based on their own intelligence assessments and internal stance. 

\subsection{Weighted threshold functions}
\label{ssec:weighted-threshold-functions}
The coalition is assumed to have a publicly known institutional design, specifically a decision rule formalized by a social choice function $f$. In order to be able to analyze the behaviour of our signalling game with respect to such a function, we need to restrict to a class of viable rules. First, we assume a binary outcome, $1$ representing coalition retaliation and $0$ representing no retaliation. Second, we assume the choice function $f$ to be a weighted threshold rule. This means that each coalition member $i$ is assigned a weight $w_i \in \R$ dependent on the context of the attack, and that there is a threshold $\tau \in \R$, such that the aggregated vote of the coalition is $1$ if the weighted sum of the individual members votes exceeds the threshold. This weight scheme is meant to include the possibility that the attacked country in the coalition gets a vote that counts more towards the total, or that frontline countries have a larger say, or that a country has far better ISR capabilities etc. 

Formally, then, given a list of the coalition members' individual votes $v=[v_1, v_2, \dots, v_n]$, and a weight scheme $w=[w_1, \dots, w_n]$, we have 
\[
f(v;w) = 
\begin{cases}
    1, &\displaystyle\sum_{i=1}^n w_i v_i \geq \tau \\
    0, &\text{else}    
\end{cases}
\]
Most viable voting rules are formalizable in this framework, like the majority rule, consensus rule, dictatorship, qualified majority, veto players, etc. There are also some decision rules that are not possible to model in this framework, like multi-stage rules, lexicographical rules, non-monotonic rules and randomized rules. However, such rules are very unlikely to be viable for choosing a retaliation action, hence our focus on this restricted class. 

\subsection{Retaliation probabilities}
Assuming $\A$ attacks, retaliation will occur if the weighted sum of coalition members' votes meets the threshold. We denote the probability of retaliation by
\[R(\tau) = \Pr\left(\sum_{i=1}^n w_i v_i \geq \tau \mid T=1 \right)\]
This function captures the \emph{deterrence credibility} of the coalition $\C$. In other words, it is the probability that an attack by the adversary $\A$ triggers a retaliatory response. 

Similarly, if the adversary does \emph{not} attack, then retaliation is accidental, with probability 
\[F(\tau) = \Pr\left(\sum_{i=1}^n w_i v_i \geq \tau \mid T=0\right)\]
This function captures the coalition's \emph{escalation risk}, i.e., the probability that the coalition accidentally authorizes a response. This could happen if, for example, another actor than $\A$ was responsible for an attack on the coalition, or that an accident is interpreted as an attack, as explored in multiple wargames and academic studies, see \cite{Sagan_1993} and \cite{Blair_1993}.

\begin{example}
    \label{ex:binomial}
    In the simple case that each coalition member $i$ has an equal weight, say $w_i = 1$ for all $i$, as well as equal probabilities $p = p_i$ and $q=q_i$, and the threshold $\tau$ is a positive integer less than or equal to $n$, then these probabilities become binomial, with 
    \[R(\tau) = \sum_{k=\tau}^{n} \binom{n}{k} p^{k} (1-p)^{n-k},\]
    and similarly 
    \[F(\tau) = \sum_{k=\tau}^{n} \binom{n}{k} q^{k} (1-q)^{n-k}.\]
The often used majority rules and consensus rules are both of this nature. 
\end{example}

\subsection{Costs and benefits}
The preferred outcome of the adversary actor $\A$ of type $T=1$ is to attack and suffer no consequences, while for $T=0$ it is to preserve the peace. The preferred outcome for the coalition $\C$ is that no attack occurs. 

In order for an attack to be of interest, there has to be a perceived benefit for $\A$. In the situation where $\A$ attacks and no retaliation occurs, we denote this benefit by $B$. We do not assume anything specific about this, but it could be monetary benefit or some other measure of benefit that $\A$ deems desirable. If retaliation occurs, then $\A$ suffers a cost $C$, which we assume is of comparable type to the benefit $B$. Note that this is intentionally simplified. 

We can then compute the expected payoff for $\A$ attacking. The probability of retaliation by $\C$ is given as above by $R(\tau)$, hence the expected payoff for $\A$ with type $T=1$ for attacking is 
\[\E(\attack) = (1-R(\tau))B-R(\tau)C.\]

\subsection{Attack conditions}
As our model is very simple from a game-theoretical perspective, the attack conditions are quite easy to find. In our model, only $\A$ is a strategic actor, who makes a choice of an action based on some information -- the coalition $\C$ does not perform a strategic action, it is decided upon by their predetermined institutional design. 

Assuming rationality, the adversary state attacks if it deems the expected payoff $\E(\attack)$ to be greater than $0$. Rearanging the terms gives that $\A$ attacks if and only if
\[R(\tau)< \frac{B}{B+C},\]
and is otherwise successfully deterred. This fractions is the \emph{deterrence threshold} for the adversary $\A$. 

We want to note that this setup is deliberately simplified, and not realistic. The setup deliberately collapses the condition into a single probabilistic parameter. We do this in order to better isolate the contribution from the institutional design to this parameter, computing how different choices for voting schemes $f$ affects the probabilities $R(\tau)$ and $F(\tau)$. In later work we want to make the above setup more realistic, but for this first study of the interactions between deterrence credibility and institutional design, we have opted to keep the model and mathematics simple. 

\begin{remark}
    We will not use this in any degree, but these two strategies and beliefs constitute the perfect Bayesian equilibria for this simplified game. Due to the simple nature, both are \emph{separating} equilibria. They are also the only possible ones, and thus yields a full classification. 
\end{remark}

\begin{remark}
    The cost $C$ depends on the coalition $\C$'s available deterrence technology. For example, having strategic nuclear weapons significantly increases the costs of retaliation for the adversary $\A$. As discussed in the introduction, new nuclear alliances are prohibited by the Non-Proliferation Treaty \cite{NPT_1968}. But, in certain situations and contexts, it is possible to circumvent this, or alternatively use other available or future technologies as in \cite{EuropeanNuclearPlanningGroup}. We can then study two different costs -- one with and one without sufficient deterrence technologies -- and see how this affects the equilibria. Doing this reveals the natural outcome: larger costs means smaller expected payoffs, and larger intervals where deterrence is successful. 
\end{remark}

As the game is very simple, it is not only the attack conditions that are interesting. The interesting question about our setup is, in our opinion, trying to determine which social choice functions that yield the best deterrence properties for the coalition $\C$. In our setup, this is also the only part of the attack condition purely in the coltrol of $\C$. The rest of the paper focuses on exploring and analyzing this question. Our method of choice does not utilize the full strategic information available, but is based purely on the probabilities of ``true'' and ``false'' retaliations. This is to isolate the effect of the institutional design on these probabilities as a first step towards a broader understanding of the full dynamics.

\section{Good schemes and thresholds}
For a given coalition $\C$, there are countless possible choices for an institutional design $f$ that computes the aggregated decision based on coalition members' votes. As mentioned earlier, this function $f$ always introduces a tradeoff between credible deterrence and escalation risk. 

For $\C$ then, it is highly important to understand which such designs that give the best deterrence credibility while minimizing the risk of accidental retaliation. In essence, this is the classical tradeoff between the probability of \emph{true positives} and the probability of \emph{false positives}. There is a rich history from many scientific diciplines studying this exact tradeoff, see for example \cite{Green-Swets_1966}, that we can use to study this issue in our setting.

Our choice function $f$ has two possible outcomes, retaliation or no retaliation, based on two possible situations of true or false attacks. By false attack, or no attack, we here mean that there has been a false detection due to an accident, or an attack that is misattributed to the wrong adversary. This means that the function $f$ is a binary classifier -- see for example \cite{Hastie-Tibshirani-Friedman_2009}. 

The standard way of measuring the accuracy of a binary classifier is to measure the correct and incorrect assignments, \cite{Swets_1988}. For this we have a simple contingency matrix:

\begin{table}[h]
\centering
\begin{tabular}{l|ll}
    & Retaliation & No retaliation \\ \hline
    Attack & True positive & False negative       \\
    No attack & False positive & True negative
\end{tabular}
\end{table}

This gives us the \emph{true positive rate} (TPR), defined as 
\[\mathrm{TPR} = \frac{\text{true positive}}{\text{true positive} + \text{false negative}}\]
and the \emph{false positive rate} (FPR), defined as 
\[\mathrm{FPR} = \frac{\text{false positive}}{\text{false positive} + \text{true negative}}.\]
These are, respectively, the probability that retaliation occurs, and the probability of an accidental retaliation, which we have previously named $R(\tau)$ and $F(\tau)$. 

The binary classifier $f$ is dependent on the threshold $\tau$, which allows us to visualize and analyze the classifiers performance on the true positive-false positive tradeoff by ROC curves. 

\subsection{ROC curves}
Receiver Operating Characteristic (ROC) curves were first used after the second world war in order to strengthen the predicting power of radar signals, see \cite{Fawcett_2006} for an introduction. Since then they have found a wide use throughout science, mainly in medicine and epidemiology \cite{Hanley-McNeil_1982}, and more recently as a tool for evaluating classification algorithms in machine learning and artificial intelligence \cite{Bishop_2006, Bradley_1997}. 

The ROC space is defined by having TPR and FPR as its axes, and as these are rates from $0$ to $1$ the space is equivalent to the unit square. Any point in this space represents a tradeoff between true positives and false positives. The point $(0,1)$ is the \emph{perfect classifier}, in the sense that it correctly identifies $100\%$, with $0\%$ errors. The diagonal line $(x,x)$ represents all random classifiers. Hence, points above the diagonal are better than random guessing, and points below are worse. 

By plotting $(F(\tau), R(\tau))$ for all different thresholds $\tau$, one gets a curve through ROC space that visualizes the effect of the threshold on the true positive-false positive tradeoff. The point on this curve that is closest to the perfect classifier $(0,1)$ is the optimal threshold for the given choice function. 

By making ROC curves for a variety of choice functions, we can then identify institutional designs for the coalition $\C$ that performs well in a variety of settings. One method for doing this is comparing ROC curves on their area under the curve (AUC) -- a better institutional design will have a larger AUC. The AUC measures the probability that the given voting scheme ranks a true positive above a false positive \cite{Hanley-McNeil_1982}, or in our case that it is able to clearly separate attacks from a known actor from ``false alarms''. Hence, a large AUC -- one that is close to $1$ -- will be more likely to produce large weighted sums when there is an actual attack, and low weighted sums when there is not. 

After finding voting schemes that perform well across all thresholds, we study which threshold that gives the best tradeoff between credible deterrence and escalation risk. We do this by comparing all thresholds on Youden's $J$ statistic, introduced originally in \cite{Youden_1950}. It is defined by $J(\tau) = R(\tau)-F(\tau)$, and is widely used in modern analyses, see for example \cite{Fluss_2005}. The threshold with the largest $J$-number has maximal vertical distance from the diagonal in the ROC space, and is interpreted as having the overall best detection quality. 

We do this in detail for an example of a small coalition. The small number of members allows us to compute the ROC curves empirically for a variety of choice functions $f$. One can also easily do this for large coalitions by using normal approximation to the probabilities, and use binormal ROC curves \cite{Metz_1978}. We have chosen not to focus on large coalitions in this paper, as results such as Condorcet's Jury Theorem \cite{Condorcet_1785, List-Goodin_2001} make these binary classification problems fairly trivial. 

\subsection{The needed computations}
The votes $v_i \in \{0,1\}$ are assumed to be independent, with probabilities $p_i$ and $q_i$ for the two types $T=\{0,1\}$. Hence, the probability of a given vote-vector $v$ is a product of independent Bernoulli variables, 
\[\Pr(v \mid T=1) = \prod_{i=1}^n p_i^{v_i}(1-p_i)^{1-v_i},\]
for the aggressive type $T=1$, and similarly 
\[\Pr(v \mid T=0) = \prod_{i=1}^n q_i^{v_i}(1-q_i)^{1-v_i}\]
for the benign type $T=0$. 

Given a weight scheme $w=[w_1, \dots, w_n]$, we have for each vote-vector $v$ a weighted sum $S=\sum_{i=1}^n w_i v_i$. There are $2^n$ possible vote-vectors, hence the distribution of this sum can be computed by aggregating all of the probabilities for all vote-vectors -- with probabilities as above -- producing two functions,
\[\Pr(S = s \mid T=1),\quad \Pr(S = s \mid T=0).\]
Computing these probabilities then, allows us to compute the points $(F(\tau), R(\tau))$ for all values of $\tau$ by computing the above probabilities for all values $s$. 

The thresholds themselves are limited to existing in the interval $[0, \sum_{i=1}^n w_i]$, as these thresholds are the only ones making an effect on the choice function. The computations in this paper were done with a simple python script, which can be downloaded from the authors GitHub page \cite{GitHub}. The program samples $\tau$ linearly $40-50$ times, which is sufficient for our simple case. 

\subsection{The example}
\label{ssec:the-example}
Our small coalition $\C$ consist of four members. We will not go into specifics or create a realistic coalition here, as this paper only serves to demonstrate the effect of the institutional design on the deterrence mechanism, and not imply policy or strategy. This also allow us to test hypothetical alterations of the choice functions due to different contexts the coalition exists in; specifically we consider choice function weight vectors $w$ based on hypothetical geographical, technological and geopolitical contexts. 

In order to compare choice functions based on their performance, we plot ROC-curves for several natural weight schemes at the same time. We have chosen seven different weight schemes for our coalition, which we will use for all of the probability variations. We have normalized them so that they all have the same total sum of $4$, as well as put the weights in decreasing order for simplicity. 
\begin{enumerate}
    \item $w = [1,1,1,1]$ -- the \emph{unbiased} scheme. This treats all coalition members equally. Varying the threshold $\tau$ then gives the standard voting schemes described in \cref{sec:social-choice-functions}. 
    \item $w = [4,0,0,0]$ -- the \emph{dictator} scheme. This is essentially the shared capability structure employed by NATO today, where the country owning the deterrence technology also has full control over its use.  
    \item $w = [2.5, 0.5, 0.5, 0.5]$ -- the \emph{veto} scheme. The scheme allows a single member to have full veto rights for most thresholds. 
    \item $w = [1.2, 1.1, 0.9, 0.8]$ -- the \emph{technology} scheme. This weighs members according to bias coming from asymmetric qualities, like techonology, competencies, or ISR capabilities. 
    \item $w = [1.6, 0.8, 0.8, 0.8]$ -- the \emph{frontline} scheme. This scheme weighs, for example, a coalition member that is a border-state with the adversary higher than the other members. 
    \item $w = [1.6, 1.2, 0.8, 0.4]$ -- the \emph{geographical} scheme. The scheme weighs members according to their distance from the aggressor, and hence by the probability of being affected. 
    \item $w = [1.3, 1.3, 0.7, 0.7]$ -- the \emph{two-bloc} scheme. This represents a coalition where there are two main members forming a block, and two secondary members forming a second block. 
\end{enumerate}

\begin{figure}[h!]
    \centering
    \includegraphics[width=\columnwidth]{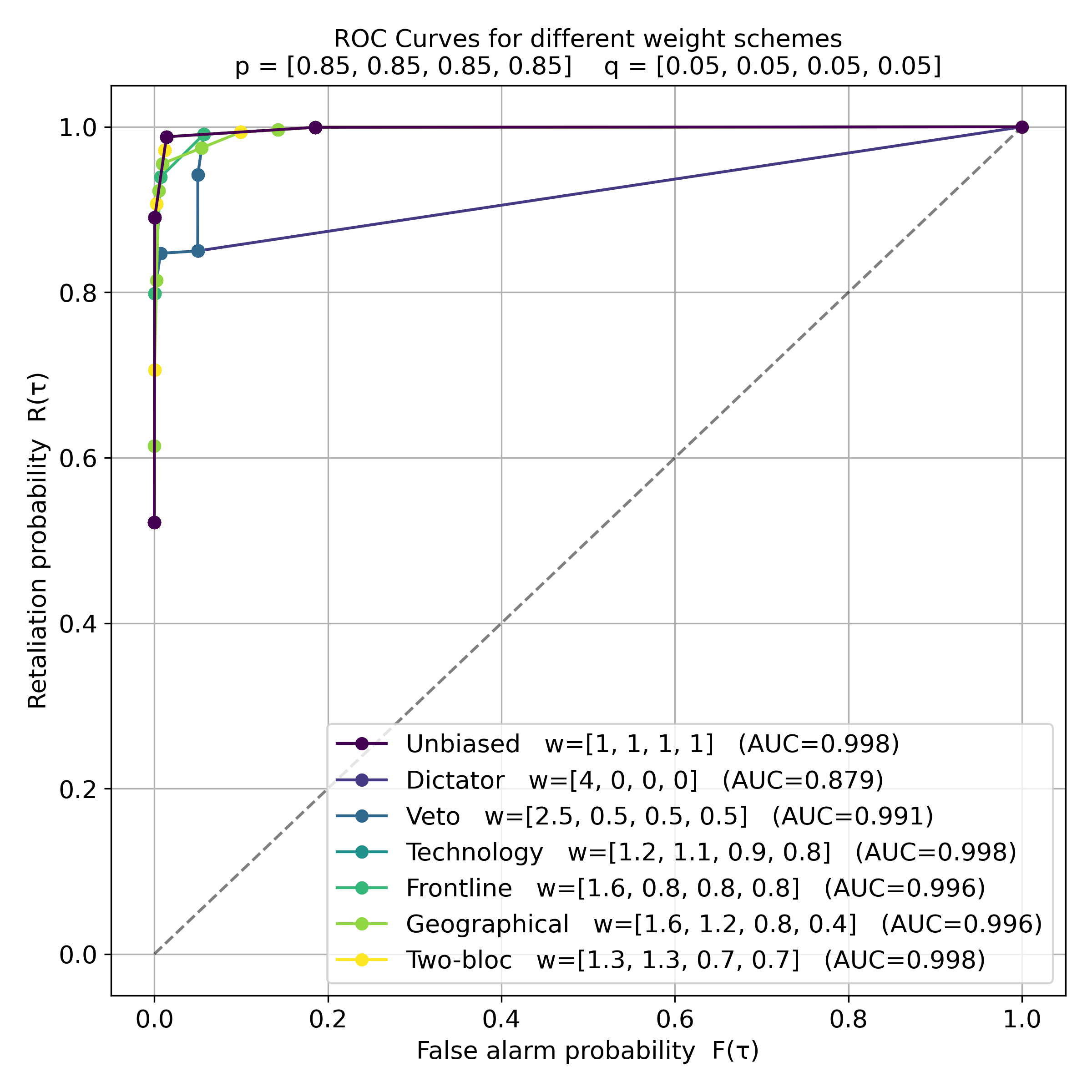}
    \caption{ROC curves for high accuracy and resolve.}
    \label{fig:high-high}
\end{figure}

To have a base line, we look at the case where all coalition members have the same probability of voting in favor of retaliation given the type of $\A$. \Cref{fig:high-high} shows the ROC curves for the above weight schemes when both the resolve and the accuracy for correct identification are \emph{high} -- represented by $p = [0.85, 0.85, 0.85, 0.85]$ and $q = [0.05, 0.05, 0.05, 0.05]$. 

As is perhaps not too surprising, all choice functions perform fairly well when the probabilities are very simple to distinguish. All AUCs are close to $1$, meaning that the probability of ranking an attack above a false alarm is almost guaranteed for all of the schemes. The unbiased scheme, the technology scheme and the two-bloc scheme performs the best, with an AUC of $0.998$, while the dictator scheme performs the worst, with an AUC of $0.879$, which is still not too bad. 

\begin{figure}[h!]
    \centering
    \includegraphics[width=\columnwidth]{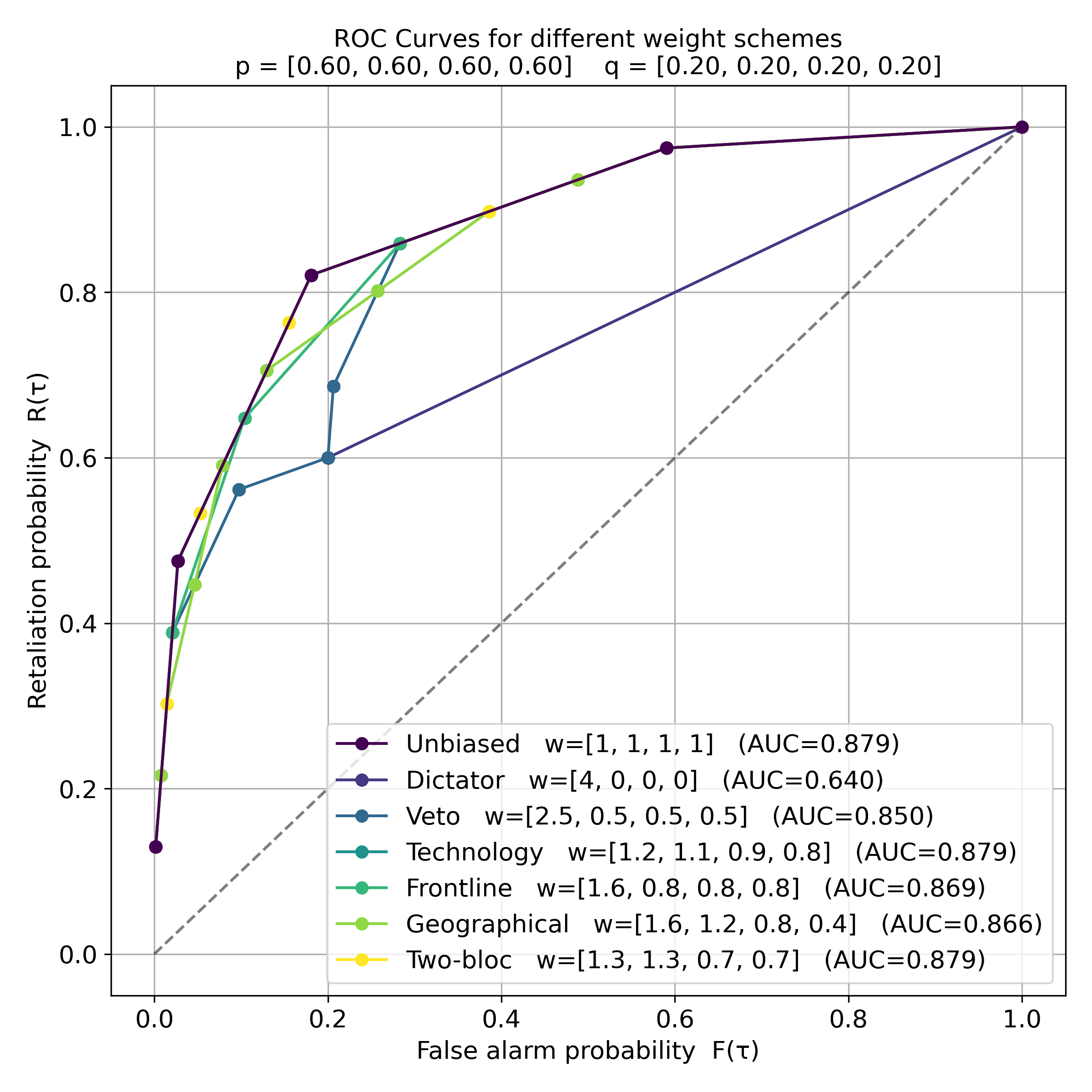}
    \caption{ROC curves for low accuracy and resolve.}
    \label{fig:low-low}
\end{figure}

Similarly, we consider in \Cref{fig:low-low} the case when the resolve and accuracy is fairly \emph{low} -- represented by $p = [0.60, 0.60, 0.60, 0.60]$ and $q = [0.20, 0.20, 0.20, 0.20]$.

As each probability is still equal for all members, the shape of the ROC curves are identical to the ones above, only pushed in closer to the diagonal. Most of the voting schemes still perform fairly well given the circumstances. 

\begin{figure}[h!]
    \centering
    \includegraphics[width=\columnwidth]{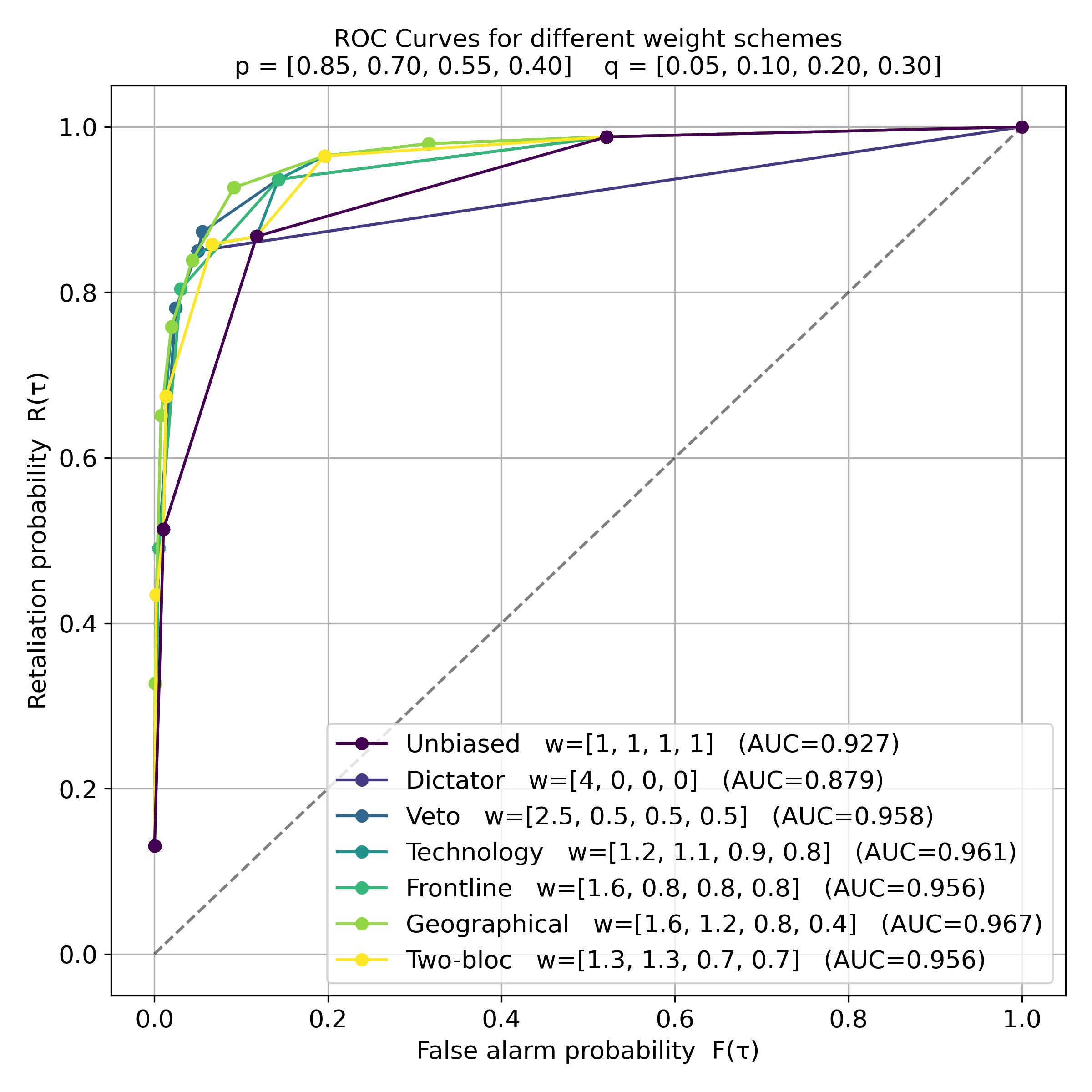}
    \caption{ROC curves for decreasing $p$.}
    \label{fig:strong-decreasing}
\end{figure}

To differentiate a bit more, we consider the following decreasing probability-vector $p = [0.85, 0.70, 0.55, 0.40]$, with the corresponding increasing vector $q = [0.05, 0.10, 0.20, 0.30]$. These probabilities align with the ranking of the members in all the weight schemes that are not the unbiased one, meaning that they ``match'' the hypothetical scenario. We then get the ROC curves in \Cref{fig:strong-decreasing}. 

Here we see that the voting schemes have different performances, and that the unbiased scheme is no longer the top performer. The best performing scheme is the one that best utilizes the distribution of the probabilities -- the geographical scheme. The dictator scheme still performs the worst. The technology scheme also performs better than the unbiased one, due to the better utilization of the probability environment. 

\begin{figure}[h!]
    \centering
    \includegraphics[width=\columnwidth]{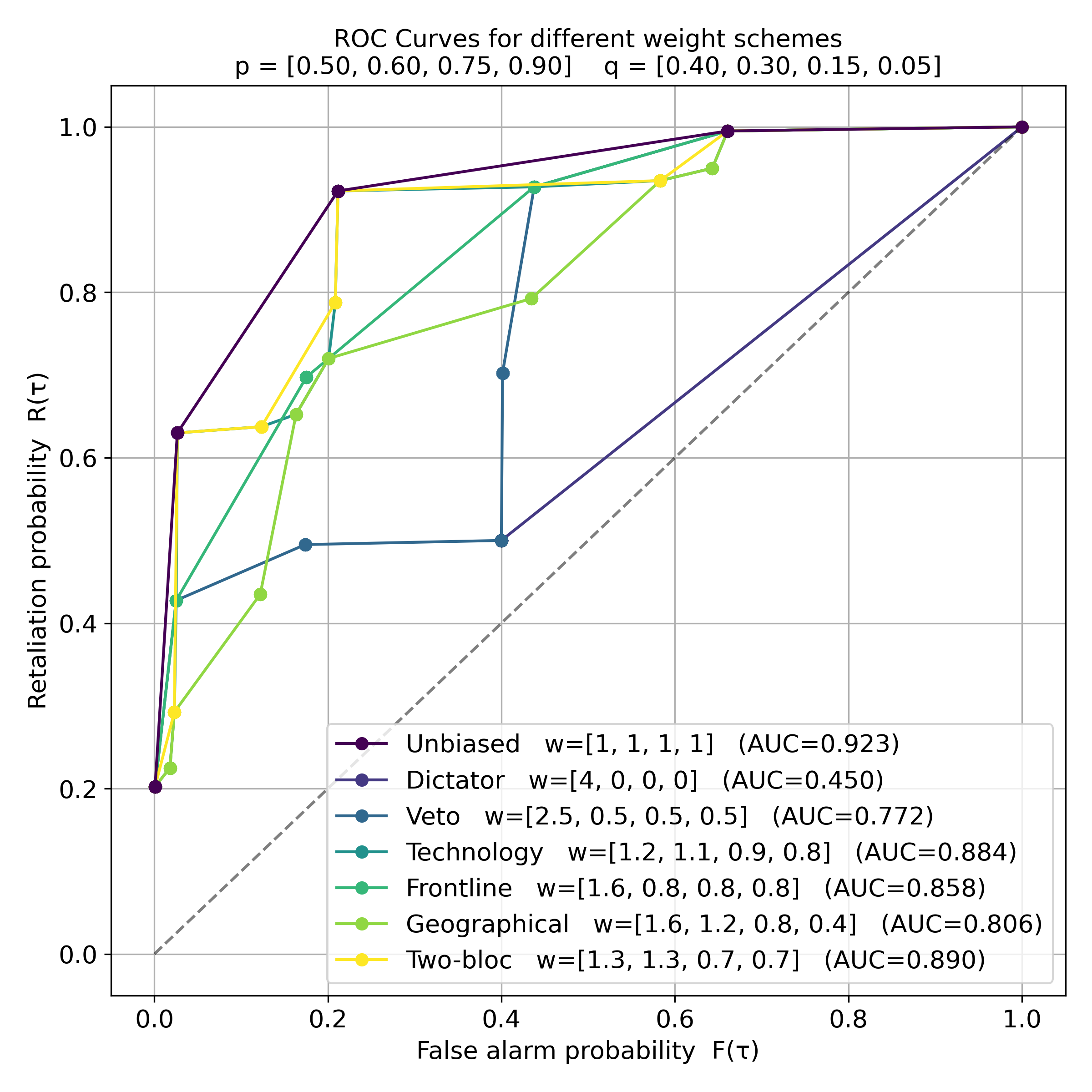}
    \caption{ROC curves for increasing $p$.}
    \label{fig:strong-increasing}
\end{figure}

We also consider the case when the probabilities are not oriented in the same direction as the weights, using the ``misaligned'' vectors $p = [0.50, 0.60, 0.75, 0.90]$ and $q = [0.40, 0.30, 0.15, 0.05]$, see \Cref{fig:strong-increasing}. 

This weight-probability combination exaggerates differences in voting schemes, which is also noticeable in the shape of the ROC curves -- they are considerably more varied. The best performer is again the unbiased scheme. The worst performer is still the dictator scheme, which here has an AUC of under $0.5$, meaning it is essentially worse than random guessing when trying to distinguish attacks from false alarms. 

\subsection{Scheme performance}
As one can see, some of the weight schemes are performing consistently better than others. To further study this we evaluate the weight schemes on a large variety of probability vectors, listed below, and study the statistics of their AUCs. These form varying \emph{information environments} that the coalition operates in. We use the probabilities as listed in \cref{table:p-q-vectors}, which are intentionally stylized for analytical purposes. 

\begin{table}[h!]
\centering
\begin{tabular}{c c}
\hline
\textbf{$p$ vector} & \textbf{$q$ vector} \\
\hline
$[0.85, 0.85, 0.85, 0.85]$ & $[0.05, 0.05, 0.05, 0.05]$ \\
$[0.75, 0.75, 0.75, 0.75]$ & $[0.10, 0.10, 0.10, 0.10]$ \\
$[0.60, 0.60, 0.60, 0.60]$ & $[0.20, 0.20, 0.20, 0.20]$ \\
$[0.90, 0.80, 0.70, 0.60]$ & $[0.05, 0.10, 0.15, 0.20]$ \\
$[0.85, 0.70, 0.55, 0.40]$ & $[0.05, 0.10, 0.20, 0.30]$ \\
$[0.60, 0.65, 0.70, 0.75]$ & $[0.25, 0.20, 0.15, 0.10]$ \\
$[0.55, 0.65, 0.75, 0.85]$ & $[0.35, 0.25, 0.15, 0.05]$ \\
$[0.50, 0.60, 0.75, 0.90]$ & $[0.40, 0.30, 0.15, 0.05]$ \\
$[0.55, 0.55, 0.55, 0.55]$ & $[0.45, 0.45, 0.45, 0.45]$ \\
$[0.85, 0.55, 0.80, 0.50]$ & $[0.05, 0.25, 0.10, 0.30]$ \\
$[0.90, 0.88, 0.60, 0.50]$ & $[0.05, 0.06, 0.25, 0.30]$ \\
$[0.72, 0.58, 0.81, 0.63]$ & $[0.11, 0.22, 0.07, 0.19]$ \\
$[0.78, 0.82, 0.79, 0.40]$ & $[0.12, 0.10, 0.11, 0.35]$ \\
$[0.88, 0.55, 0.70, 0.83]$ & $[0.06, 0.28, 0.18, 0.10]$ \\
\hline
\end{tabular}
\caption{Varied collection of $p, q$ pairs.}
\label{table:p-q-vectors}
\end{table}

Computing AUC values for all of these, we obtain the statistics for each weight scheme as in \cref{table:AUC-statistics}. 

\begin{table}[h!]
\centering
\begin{tabular}{lccc}
\hline
\textbf{Scheme} & \textbf{Mean} & \textbf{Min} & \textbf{Max} \\
\hline
Unbiased     & 0.932 & 0.606 & 0.998 \\
Dictator     & 0.719 & 0.426 & 0.903 \\
Veto         & 0.902 & 0.594 & 0.991 \\
Technology   & 0.931 & 0.606 & 0.998 \\
Frontline    & 0.922 & 0.602 & 0.996 \\
Geographical & 0.914 & 0.600 & 0.988 \\
Two-bloc     & 0.929 & 0.606 & 0.998 \\
\hline
\end{tabular}
\caption{Mean, minimum and maximum AUC values, computed from the $p,q$ pairs in \cref{table:p-q-vectors}.}
\label{table:AUC-statistics}
\end{table}

The best performer is the unbiased weight scheme. It has the highest mean AUC, as well as the highest minimum AUC, and the highest maximum AUC. It is, however, followed closely by the technology scheme and the two-bloc scheme, which shares its minimum and maximum. The dictator scheme performs considerably worse than the others, with the lowest mean, the lowest minimum and the lowest maximum AUC. 

\begin{figure}[h!]
    \centering
    \includegraphics[width=\columnwidth]{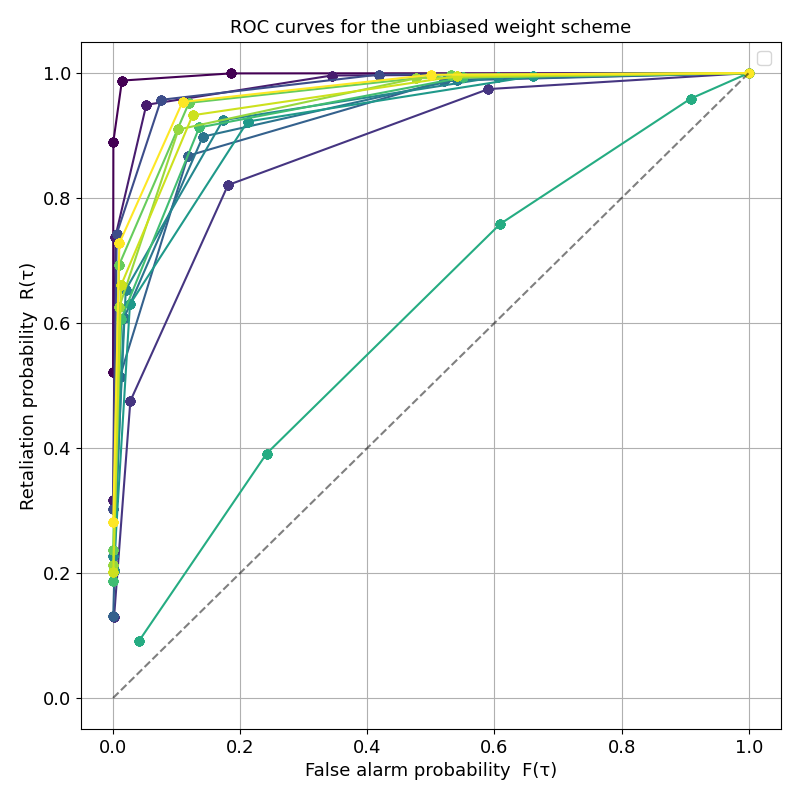}
    \caption{ROC curves for the unbiased scheme on the $p,q$ pairs in \cref{table:p-q-vectors}.}
    \label{fig:unbiased}
\end{figure}

To visualize the performance of the unbiased scheme, we map all of its ROC curves against all the probability vectors in \cref{table:p-q-vectors} in \cref{fig:unbiased}. 

We see that its shape is consistent, and that its performance is good for most of the probability vectors. The curve close to the diagonal is associated to the vectors $p = [0.55, 0.55, 0.55, 0.55]$ and $q= [0.45, 0.45, 0.45, 0.45]$, which is close to being a random guess for the type of the adversary $\A$. Note that, if $p_i>0.5$ and $q_i<0.5$, Condorcet's Jury theorem (\cite{Condorcet_1785, List-Goodin_2001}) assures that the curve approaches the perfect classifier $(0,1)$ as the sixe of the coalition increases. 

\begin{figure}[h!]
    \centering
    \includegraphics[width=\columnwidth]{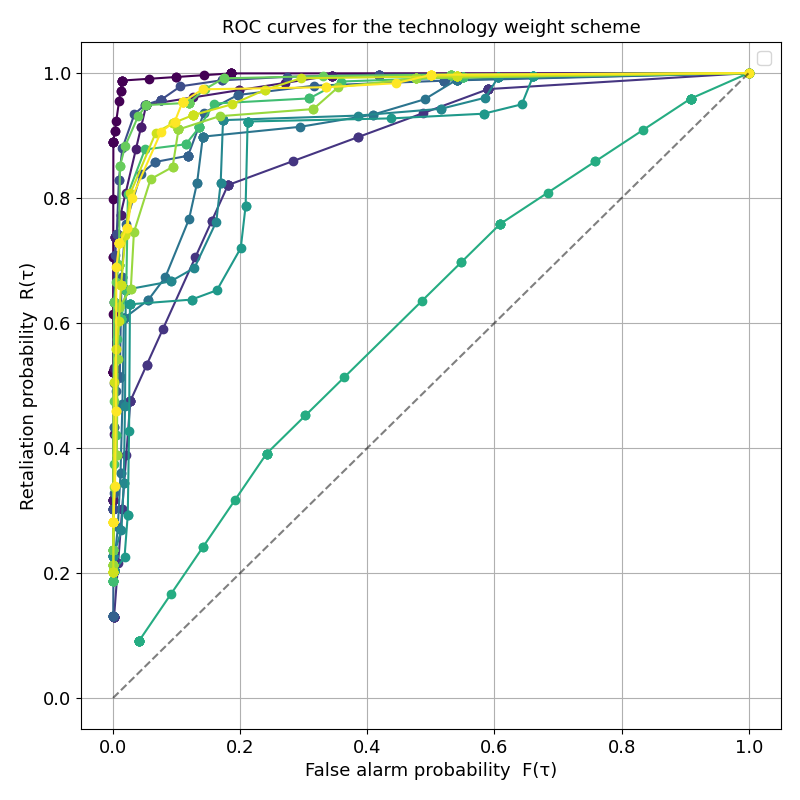}
    \caption{ROC curves for the technology scheme on the $p,q$ pairs in \cref{table:p-q-vectors}.}
    \label{fig:technology}
\end{figure}

Similarly, for the technology scheme we get the collection of ROC curves presented in \cref{fig:technology}. 

Its performance is statistically similar to the unbiased scheme, but we see that it is much more sensitive to threshold variations, and that the resulting ROC curves are then more varied in shape. These alternative shapes come from probability vectors that are misaligned with the weight gradient -- which is purposefully not representing the schemes power in utilizing individual differences. This would be like a member getting higher voting power for having better ISR technology, but then being consistently unable to utilize this technology to its full potential, while other members with worse technology manage it better. 

%Finally, \cref{fig:dictator} shows the same data for the worst performer -- the dictator scheme. 

%\begin{figure}[h!]
%    \centering
%    \includegraphics[width=\columnwidth]{images/dictator_no-label.png}
%    \caption{ROC curves for the dictator scheme on the $p,q$ pairs in \cref{table:p-q-vectors}.}
%    \label{fig:dictator}
%\end{figure}

%Its ROC curves are all shaped like lines, as it is only possible for the scheme to output a $1$ when the threshold is lower than the voting power of the dictator -- all other thresholds give the same point $(1,1)$. Note that the dictator scheme can perform well compared to the other schemes, in situations where the dictator member is vastly superior in information.

\subsection{Threshold performance}
The previous section ranked each weight scheme on their performance across \emph{all} thresholds. This is, however, not the final piece of the puzzle, a we still need to connect performance to specific thresholds $\tau$. In essence we want to find the weight-threshold \emph{combination} that gives the best performing single point in a our varied information environments -- not the entire curve. 

To do this we compute Youden's $J$-statistic for all schemes on all of our earlier information environments, and study the statistics of the schemes performances. Youden's $J$-statistic is, as previously mentioned, defined by $J(\tau) = R(\tau) - F(\tau)$, and measures the greatest vertical distance from the random classifier diagonal line. The threshold that yields the greatest $J$-number is the \emph{optimal} threshold, and will be denoted by $\tau^*$. 

Just as for the AUC values, we list in \cref{table:J-statistics} the mean, minimum and maximum $J$-numbers for all of the voting schemes, computed for all of the probability vectors in \cref{table:p-q-vectors}. We also list the mean optimal threshold $\tau^*$, computed for the same data. The $\tau^*$ value below should really be discretized by rounding up to the nearest whole integer, but we have decided to keep them as is for comparison and insight.  

\begin{table}[h!]
\centering
\begin{tabular}{lcccc}
\hline
\textbf{Scheme} 
& \textbf{Mean} 
& \textbf{Min} 
& \textbf{Max} 
& \textbf{Mean $\tau^*$} \\
\hline
Unbiased      
& 0.807 
& 0.150 
& 0.974 
& 1.17 \\

Dictator      
& 0.566 
& 0.100 
& 0.850 
& 0.10 \\

Veto          
& 0.705 
& 0.125 
& 0.934 
& 0.84 \\

Technology    
& 0.802 
& 0.150 
& 0.974 
& 1.44 \\

Frontline     
& 0.757 
& 0.149 
& 0.934 
& 1.19 \\

Geographical  
& 0.760 
& 0.149 
& 0.946 
& 1.48 \\

Two-bloc      
& 0.807 
& 0.150 
& 0.974 
& 1.47 \\
\hline
\end{tabular}
\caption{The mean, minimum and maximum Youden's $J$ statistic, and mean optimal threshold $\tau^*$ for all voting schemes, computed across the $p,q$ pairs in \cref{table:p-q-vectors}.}
\label{table:J-statistics}
\end{table}

Here we see that the unbiased scheme and the two-bloc scheme performs the best overall, followed closely again by the technology scheme. The two-bloc scheme has a higher mean optimal threshold $\tau^*$ compared to the unbiased scheme, but in effect, due to the discrete small coalition, the best threshold is $\tau^* = 2$ for both. Only the dictator and veto schemes has an optimal threshold of $\tau^* = 1$. 

As the unbiased scheme does not make any assumptions on the coalition's internal structure -- the partitioning into blocs -- as well as having better AUC statistics, it is the preferred voting scheme for the four-member deterrence coalition $\C$ -- this is precisely the \emph{majority rule}. 

Intuitively, the majority rule benefits from its essential feature: symmetry. When a member's individual errors are independent and not systematically biased, the unbiased scheme's equal weighting minimizes sensitivity to these miscalibrations, making it perform best on average. In the case where realistic biases are thoroughly mapped and utilized, this symmetry reduces the effectiveness and advantage of biased weights, but in an uncertain world where all things can happen, symmetry seems to be the best approach.  

\begin{remark}
    If the coalition  is highly certain and confident of their members technological capabilities and its accuracy in a plethora of real-world scenarios, then the technology scheme might, as we saw in \cref{ssec:the-example}, outperform it -- still with optimal threshold $\tau^* = 2$. However, given any possibility for a scenario where the probabilities are misaligned with the weights, the unbiased scheme will outperform the technology scheme on average. 
\end{remark}

\section{Conclusion}
Now that we have landed on a well-performing social choice function for our coalition $\C$, we can also understand the associated deterrence behavior for the signaling game-inspired model in \cref{sec:the-model}. 

If we omit the unlikely ``outlier'' scenario that the probability vectors are almost random, $p=[0.55,0.55,0.55,0.55]$ and $q=[0.45, 0.45, 0.45,0.45]$, we obtain an average retaliation probability 
\begin{align*}
    \overline{R(\tau^*)}
    &= \frac{1}{13}\sum_{k=1}^{13}\Pr\left(\sum_{i=1}^4 v_i \geq 2 \mid T=1, p = p_k, q=q_k\right) \\
    &\approx 0.92
\end{align*}
where $p$ and $q$ ranges over the remaining $13$ probability vectors in \cref{table:p-q-vectors}. Similarly, 
\begin{align*}
    \overline{F(\tau^*)}
    &= \frac{1}{13}\sum_{k=1}^{13}\Pr\left(\sum_{i=1}^4 v_i \geq 2 \mid T=0, p = p_k, q=q_k\right) \\
    &\approx 0.12.
\end{align*}
On average, across a plethora of scenarios, the adversary state $\A$ attacks when the perceived benefit and costs satisfy $0.92 \leq \frac{B}{B+C}$, which happens if the benefit is about $11.5$ times greater than the costs. Hence, for four-member coalitions with sufficient deterrence technologies, like strategic nuclear warheads, deterrence will be quite credible in most situations when using the majority rule as their institutional design.

In summary, we have presented a simple model for understanding an adversary's attack conditions, dependent on how a deterrence coalition aggregates its members votes toward a shared decision on retaliation. This, we argued, introduced a strategic tradeoff between true and false positives, or \emph{deterrence credibility} and \emph{escalation risk}. This tradeoff was used to further frame deterrence as a binary classification problem, which allowed us to use standards techniques from medical analysis and signal analysis to explore how different such institutional designs respond to different information environments, and how it affects the adversary's beliefs about retaliation.

\printbibliography{}  
\end{document}